\documentclass[12pt]{amsart}
\usepackage{amscd, amssymb, amsthm}
\numberwithin{equation}{section}

\newtheorem{pr}[equation]{Proposition}

\newtheorem{lm}[equation]{Lemma}
\newtheorem{cor}[equation]{Corollary}
\newtheorem{tm}[equation]{Theorem}

\newcommand{\proj}{{\mathbf P}}

\newcommand{\barr}{\overline}
\newcommand{\rarr}{\rightarrow}
\newcommand{\oh}{{\mathcal{O}}}
\newcommand{\CC}{{\mathcal{C}}}
\newcommand{\EE}{{\mathcal{E}}}
\newcommand{\RR}{{\mathcal{R}}}
\newcommand{\A}{{\mathbb{A}}}
\newcommand{\F}{{\mathbb{F}}}
\newcommand{\C}{{\mathbb{C}}}
\newcommand{\Q}{{\mathbb{Q}}}
\newcommand{\Z}{{\mathbb{Z}}}

\newcommand{\PGL}{{\mathbf{PGL}}}

\newcommand{\sumo}{\oplus}
\newcommand{\M}{\barr{M}}
\newcommand{\Mon}{\M_{0,n}}
\newcommand{\bpf}{\begin{proof}}
\newcommand{\epf}{\end{proof}}
\newcommand{\Mgn}{\M_{g,n}}
\newcommand{\smcup}{\mathbin{\text{\scriptsize$\cup$}}}

\newcommand{\MonX}{\M_{0,n}(X,\beta)}
\newcommand{\MgnX}{\M_{g,n}(X,\beta)}
\newcommand{\MgoX}{\M_{g,0}(X,\beta)}

\newcommand{\MonH}{\M_{0,n}(H,\beta)}

\newcommand{\mooH}{M_{0,0}(H,\beta)}

\newcommand{\Sym}{{\mbox{Sym}}}
\newcommand{\Pic}{{\mbox{Pic}}}
\newcommand{\Ext}{{\mbox{Ext}}}
\newcommand{\Spec}{{\mbox{Spec}}}
\newcommand{\Hom}{{\mbox{Hom}}}

\newcommand{\smprod}{\prod}

\newcommand{\XX}{{\mathcal Z}}

\def\fra#1#2{{{#1}\atop{#2}}}

\def\feyn#1#2#3#4{\displaystyle\Bigl(
\fra{\displaystyle\hfill #1}{\displaystyle\hfill #2}{\rangle\hskip-3pt}
\frac{\ \ }{\ \ }{\hskip-3pt\langle}
\fra{\displaystyle #4\hfill}{\displaystyle #3\hfill}\Bigr)}

\begin{document}

\title{Enumerative Geometry of Hyperelliptic Plane Curves}
\author{Tom Graber}

\maketitle

\setcounter{section}{-1}

\section{Introduction}

In recent years there has been a tremendous amount of progress
on classical problems in enumerative geometry.  This has largely
been a result of new ideas and motivation for these problems
coming from theoretical physics.  In particular, the theory
of Gromov-Witten invariants has 
provided powerful tools for counting curves satisfying 
incidence conditions.

This theory has been most successful in dealing with questions
about rational curves.  This is partly because
it is much more common for the genus 0 invariants to
correspond to enumerative problems.  In addition, it
is much easier to compute these invariants.
This is due mainly to the existence
of the WDVV equations. 
There has been success in extending
the techniques used to derive these equations to
find recursions satisfied by the invariants in genus 1 and 2.
(\cite{ezra},\cite{rp}).  In some situations these higher
genus invariants also correspond to classical enumerative problems.
Thus, the theory gives new methods to solve these problems.

We will explore a different approach to using Gromov-Witten
theory to solve enumerative problems involving higher genus
curves.  Rather than generalizing the methods that
succeed in genus 0, we try to reduce questions in higher genus
to questions about rational curves.
We utilize the well-developed theory of genus 0 
Gromov-Witten invariants to solve enumerative
problems involving hyperelliptic curves in $\proj^2$.
Our main enumerative result is the construction of 
a recursive algorithm which
counts the number of hyperelliptic plane curves of
degree $d$ and genus $g$ passing through $3d+1$ general
points.

The basic facts about stable maps and Gromov-Witten theory are reviewed 
in Section \ref{sectiongw}.  The main results we need are the
WDVV equations which allow us to compute, and an explicit
representation of the virtual fundamental class as
a Chern class.

In Section \ref{mysterysection} we introduce the main idea of the paper. 
Thinking of a map from a hyperelliptic curve as a family of 
maps from pairs of points parametrized by $\proj^{1}$ gives us
a natural correspondence between hyperelliptic plane curves
and rational curves in $H=H(2,\proj^{2})$, the Hilbert scheme
of two points in the plane.  Furthermore, the condition that the
hyperelliptic curve meet a point $p$ in $\proj^{2}$ is equivalent
to the condition that the associated rational curve in $H$
meets the cycle 
\[
\Gamma (p) = \{\mbox{subschemes incident to {\it p}} \}.
\]
Since genus 0 Gromov-Witten invariants naively correspond to 
counting rational curves incident to cycles, it is reasonable to
hope that we could understand our enumerative problem in terms
of Gromov-Witten theory on $H$.
We study the geometry of this Hilbert scheme, and identify
which curve classes arise via this correspondence.  

These results enable us to identify exactly how the Gromov-Witten
invariants differ from the solution to the enumerative problem
in which we are interested.  There are extraneous components
of the space of maps which contribute to the invariants.  In
order to relate the Gromov-Witten theory to enumerative geometry, it
is necessary to identify
the contributions from these components.  This
is carried out in section \ref{deform}. 
By studying the deformation theory on these components,
we can relate their contributions to a much
simpler problem about curves on a blow up of the plane
where it is trivial to determine the solution.
The conclusion is a simple formula relating the
enumerative numbers to the Gromov-Witten invariants.
If we set $E(d,g)$ to be the number of degree $d$ genus $g$
hyperelliptic curves passing through $3d+1$ general points
in the plane, and set $I(d,g)$ to be the corresponding
Gromov-Witten invariant, then
our main result is that
$$ I(d,g) = \sum_{h\geq g} \binom{2h + 2}{h-g} E(d,h).$$
This relationship can be inverted to solve the enumerative problem
in terms of Gromov-Witten invariants.  We actually get a similar
formula for the more exotic enumerative problem of counting
hyperelliptic plane curves passing through fixed points, with
certain pairs of the points required to be hyperelliptically conjugate.
This condition can be used to recover genus 0 and 1 Severi degrees
from this formalism.
We also
use our calculations to recover a formula of Abramovich and Bertram
\cite{ab} about genus 0 Severi degrees on the Hirzebruch surfaces
$\F_0$ and $\F_2$.  

In Section \ref{calc} we show how
we can use the formal properties of the invariants to give 
a recursive algorithm computing the $I(d,g)$.  In this way we are 
able to effectively compute the solution to our original problem.
In addition, we use our calculations to give an explicit presentation
of the small quantum cohomology ring of $H$.

This bulk of the work presented in this paper 
was carried out during the author's stay at the Mittag-Leffler Institute.
We would like to thank the organizers of the special year in quantum
cohomology for providing a wonderful atmosphere for research.  It 
is also a pleasure to thank W. Fulton 
for making possible the author's visit to the
University of Chicago, where this
paper was completed.  
We are grateful to P. Belorousski, C. Faber, B. Fantechi, and 
A. Kresch for valuable conversations about this work.  Our main debt is
to R. Pandharipande who suggested this problem and was extremely helpful
throughout our work on it. 
This research was supported at different times 
by an NSF graduate fellowship and
a Sloan dissertation year fellowship.

\section{Gromov-Witten Invariants}
\label{sectiongw}

\subsection{Stable Maps}\label{stable maps}
The Gromov-Witten invariants of a smooth projective variety $X$  
are defined as integrals over the space of stable maps, $\MgnX $.
We refer
the reader to \cite{fp} for a construction of this space, and a
careful discussion of its properties.  We should note that while
in that reference the emphasis is on the coarse moduli scheme, we
will take all our moduli spaces to be in the category of
Deligne-Mumford stacks.  We will briefly review
some of the fundamental structures on this space to fix our notation.

We denote by $M_{g,n}(X,\beta)$ 
the open substack of $\MgnX$ which parametrizes
maps from smooth irreducible curves.  While
we think of $\MgnX$ as a compactification of this space, it
is entirely possible for a component of the compact space
to parametrize only maps from singular curves.  In other words, 
$M_{g,n}(X,\beta)$ need not be dense in $\MgnX$.

These spaces come equipped with several natural morphisms.  
For each of the $n$ marked points, there is an
evaluation map, 
$$\rho_i:\MgnX \rarr X,$$ 
which takes the moduli point
$(C;x_1, \ldots, x_n;f)$ to the point $f(x_i)$. 
There is the map 
$$ \pi_{n+1} : \M_{g,n+1}(X,\beta) \rarr \MgnX $$
which simply forgets the last marked point (and stabilizes the curve
if necessary.)  This realizes $\M_{g,n+1}(X,\beta)$ as the universal
curve over $\MgnX$.  The universal morphism is given by the
$(n+1)$st evaluation map.  
$$
\begin{CD}
 \M_{g,n+1}(X,\beta) @>{\rho_{n+1}}>>  X\\
@V{\pi_{n+1}}VV \\
\MgnX \\
\end{CD}
$$

 In fact, given any subset of $\{1, \ldots , n\}$ there
is the analogous map which forgets all marked points in that
subset.  We will only have use for this construction when the
subset is all of the points, and we will denote the
corresponding map from $\MgnX$ to $\MgoX$ by $\pi$.
Lastly, there is a natural morphism
$$\eta: \MgnX \rarr \Mgn$$ 
which simply forgets the map
(and again stabilizes the curve if necessary.)

\subsection{Deformation Theory}
Crucial to the Gromov-Witten theory is a good understanding
of the local structure of these moduli spaces.  There is a 
natural obstruction theory for stable maps which 
formally locally realizes
the moduli space as the zero locus of a collection of equations
in its Zariski tangent space.  

At a fixed moduli point, this theory consists of two vector spaces,
a tangent space and an obstruction space.  These spaces are
$\Ext^1(f^*\Omega_X \rarr \Omega_C, \oh_C)$ and 
$\Ext^2(f^*\Omega_X \rarr \Omega_C, \oh_C)$ respectively.
These varying vector spaces can also be thought of as 
global coherent sheaves on the moduli space.  
From here on, these two sheaves will
be denoted by $TM$ and $\EE$.

There are two facts about these sheaves 
that will be important for us.
The first is that
on the smooth locus of the moduli space, both are locally free.
For $TM$ this is true by definition, and it follows that $\EE$ is
locally free from the existence of a presentation for
these two sheaves as a kernel and cokernel of a two term
complex of vector bundles.

The other thing we need to know about them is that they 
fit into the following exact sequence:
\begin{equation}
\label{tanobseq}
\begin{split}
0 \rarr \Ext^0&(\Omega_C, \oh_C) \rarr H^0(f^*TX) \rarr TM \rarr \\
&\rarr \Ext^1(\Omega_C, \oh_C) \rarr H^1(f^*TX) \rarr \EE \rarr 0.
\end{split}
\end{equation}
(Here, and throughout the paper, we will often refer to sheaves
by simply naming their fibers.  So the $H^i$ should be thought of as
$R^i\pi_*$ and similarly for the $\Ext$ groups.)

To understand the geometry of this sequence, it should be observed
that the $\Ext^i(\Omega_C, \oh_C)$ are the spaces of
automorphisms and deformations of the underlying nodal curve,
and the $H^i(C,f^*TX)$ are the tangent and obstruction
spaces to the space of maps from the fixed curve $C$ to $X$.

Notice that if $H^1(f^*TX)$ vanishes, we can conclude
that $\EE$ vanishes also.  This forces the
moduli space to be smooth of the expected dimension.  
We will also want to know the following stronger fact, which
is an easy consequence of Theorem II.1.7 in \cite{kollar}.
\begin{tm}
\label{relobtm}
If $f:C\rarr X$ is a morphism from a stable pointed curve
such that $H^1(C,f^*TX)=0$,
the forgetful morphism 
$$\eta:\MgnX \rarr \Mgn$$
is smooth at $[f]$.
\end{tm}

\subsection{Virtual Fundamental Class}

The main
technical point in defining the Gromov-Witten invariants is the
construction of a natural Chow homology class in $A_*(\MgnX)$ -- the
{\em virtual fundamental class} which we denote by $[\MgnX]^{\it vir}$.  
It is pure dimensional
of the expected dimension
$$ d_X(\beta) = -K_X \cdot \beta + (\dim(X)-3)(1-g) + n.$$
Given this class,
the Gromov-Witten invariants are
defined to be multilinear maps from the cohomology of $X$ to $\Q$
given by
$$ I^X_{g,\beta}(\gamma_1 \cdots \gamma_n) = 
\int_{[\MgnX]^{\it vir}}
\rho_1^*\gamma_1 \smcup \cdots \smcup \rho_n^* \gamma_n .$$

This virtual class is constructed in \cite{lt} and in \cite{bf},
\cite{b}.  In general, the construction is quite subtle.
Fortunately, we will need to know 
relatively little.  One fact we need is
that the virtual class on the pointed map space is pulled
back from the unpointed space by the flat morphism $\pi$
which forgets the marked points.  Since it is generally easier
to think about these unpointed spaces anyway, we will do all
of our calculations there.

On the unpointed spaces, we will need to compute only the restriction
of the virtual class to the smooth locus of the moduli space.
This can be realized as the top
Chern class of the obstruction bundle $\EE$.  Essentially this
is a manifestation of the standard fact that for smooth varieties,
excess intersection classes can be described as Chern classes
of associated vector bundles.

In particular, the virtual class can naturally be thought of as
a cohomology class in this situation.  We will tend to think of
the $\gamma_i$ as homology classes, effectively reversing the 
usual formulation.  That is, we will represent the cohomology
classes $\gamma_i$ by algebraic cycles $\Gamma_i$, pull these
back to the moduli space, intersect them, and integrate the
virtual class over the resulting cycle.  As we want to work
on the unpointed space, we will actually push this cycle forward
to $\MgoX$ and integrate the appropriate Chern class over the image
(up to a possible multiplicity.)  
The formulation we ultimately want is the following.

\begin{tm}
\label{realgw}
Suppose $\Gamma_1, \ldots , \Gamma_n$ are cycles in $X$ 
representing the cohomology classes $\gamma_i$ such that
$\rho_i^{-1}(\Gamma_i)$ intersect generically 
transversally.  Then if 
$$\pi_* ([\cap_i \rho_i^{-1}(\Gamma_i)]) = A$$
where $A$ is a cycle contained in the smooth locus
of $\MgoX$,
$$I^X_{g,\beta} (\gamma_1 , \ldots , \gamma_n)=  \int _A c_{\rm top}(\EE).$$
\end{tm}

The only real content of this statement is the previously mentioned
identification
of the virtual class in this setting with the top Chern class of the
obstruction bundle.  This is Proposition 5.6 in \cite{bf}.

\subsection{Properties of Gromov-Witten Invariants}

The Gromov-Witten invariants satisfy many formal properties
which can make computing them more tractable than solving
enumerative problems in general.  

First of all, the motivating property satisfied by the invariants
is that they are invariant under deformations of $X$.  There
are several ways to express this property, but essentially
the point is that since smooth families of varieties are locally
trivial in the $C^\infty$ category, there are canonical
isomorphisms between the cohomology groups of nearby
fibers.  With respect to this identification, the Gromov-Witten
invariants are independent of the choice of fiber.

Another formal property of the invariants that we will use
is the so called {\em divisor axiom}  
$$ I^X_{g,\beta}(\gamma_1 \cdots \gamma_n) = 
(\gamma_1 \cap \beta) \cdot I^X_{g,\beta}(\gamma_2 \cdots \gamma_n).$$
Repeated use of this axiom reduces the
computation of any Gromov-Witten invariant to one which involves
no divisor classes.

The genus 0 invariants behave better in several ways than
the higher genus invariants.  Since we will work with only
the genus 0 Gromov-Witten invariants,
from now on we use $I^X_\beta(\gamma_1 \cdots \gamma_n)$ to mean 
$I^X_{0,\beta}(\gamma_1 \cdots \gamma_n)$.  We will also suppress
the $X$ from the notation except where there is possibility of confusion.

The invariants can be computed much more readily in the
genus 0 case than in general.  There are essentially two 
reasons for this.  One is that it is much more common for
the $H^1(f^*TX)$ to vanish for genus 0 maps
which means that for sufficiently nice varieties, there
are no virtual class considerations.  This is the case,
for example, for all homogeneous varieties. 

The other more general fact is that the genus 0 invariants
satisfy an important family of recursive equations which
often determine all of the invariants from a small
amount of initial data.  These relations arise by considering
the many natural morphisms from the moduli space to $\proj^1$
given by composing the map $\eta : \MonX \rarr \Mon$ with
a map from $\Mon$ to $\M_{0,4}$ given by forgetting all but 
4 of the points.  As $\M_{0,4} \cong \proj^1$, different points
on it are rationally equivalent.  By pulling back the points
corresponding to reducible curves, one obtains linear equivalences
between divisors on $\MonX$.  The WDVV equations are then
deduced by intersecting these divisors with curves to obtain
numerical equalities.  These arguments all ultimately rest
on nice properties of the virtual fundamental classes.
We refer the reader to \cite{km} or \cite{fp} for a careful discussion
of this construction in the homogeneous case,
and to \cite{b}, \cite{lt} for the additional arguments needed
to handle the virtual considerations needed in 
general.  

The final result can be described as follows.
If we choose a homogeneous basis for the cohomology of $X$ given by
$T_0, T_1, \ldots, T_m$, then we get a relation for each diagram
of the following form

$$\feyn ijkl \ \sim\ \Bigl(
\begin{picture}(40,10)(-20,17)
\put(0,16){\line(0,1){8}}
\put(0,16){\line(3,-1){10}}
\put(0,16){\line(-3,-1){10}}
\put(0,24){\line(3,1){10}}
\put(0,24){\line(-3,1){10}}
\put(12,8){\em{k}}
\put(-17,8){\em{j}}
\put(12,26){\em{l}}
\put(-17,26){\em{i}}
\end{picture}\Bigr)$$

Typically the relation is described by first constructing a generating
function using the Gromov-Witten invariants as coefficients, and then
producing a differential equation which is satisfied by that
function.  This formalism is in general very useful, but for
our purposes we really just want to know the explicit recursive
equations that arise from identifying coefficients in this expression.
To this end, we need to choose a collection of cohomology classes on $X$,
 $\gamma_1,\ldots,\gamma_n$, and a class $\beta$ in $H_2(X,\Z)$.  The
relation that is obtained is the following:

\begin{equation}
\label{WDVV}
\begin{split}
%I_\beta(T_i\cdot T_j \cdot T_k\smcup T_l\cdot \smprod_1^n \gamma_i)
%+I_\beta(T_k\cdot T_l \cdot T_i\smcup T_j \cdot \smprod_1^n \gamma_i)+
&\sum_{{{{\beta_1 +\beta_2 =\beta}\atop{A\cup B = [n]}}\atop {e,f}}}
I_{\beta_1}(T_i\cdot T_j \cdot T_e \cdot \smprod_{a\in A}\gamma_a) g^{ef}
I_{\beta_2}(T_k\cdot T_l \cdot T_f \cdot \smprod_{b\in B}\gamma_b)
\\
&= 
%I_\beta(T_i\cdot T_l \cdot T_j\smcup T_k\cdot \smprod_1^n \gamma_i)
%+I_\beta(T_j\cdot T_k \cdot T_i\smcup T_l \cdot \smprod_1^n \gamma_i)+
\sum_{{{{\beta_1 +\beta_2 =\beta} \atop {A\cup B = [n]}} \atop {e,f}}}
I_{\beta_1}(T_i\cdot T_l \cdot T_e \cdot \smprod_{a\in A}\gamma_a) g^{ef}
I_{\beta_2}(T_j\cdot T_k \cdot T_f \cdot \smprod_{b\in B}\gamma_b).
\end{split}
\end{equation}

Because the space of degree 0 maps is so simple, the terms in the
summation with $\beta_1$ or $\beta_2$ equal to 0 have a particularly
simple form.  Summing over all the remaining indices, the contribution
of these terms on the left hand side of \ref{WDVV} is just
$$
I_\beta(T_i\cdot T_j \cdot T_k\smcup T_l\cdot \smprod_1^n \gamma_i)
+I_\beta(T_k\cdot T_l \cdot T_i\smcup T_j \cdot \smprod_1^n \gamma_i).
$$
Of course, the analogous expression gives the
$\beta=0$ terms for the right hand side of \ref{WDVV} as well.
Using these equations, as well as the divisor axiom, Kontsevich
and Manin show in \cite{km} 
that on a variety whose cohomology is generated
by divisors, all the genus zero 
Gromov-Witten invariants can be recursively
determined by knowledge of the two point invariants, $I_{\beta }(\gamma _{1}
\cdot \gamma _{2})$.

\section{The Hilbert Scheme}
\label{mysterysection}
\subsection{Geometry of the Hilbert Scheme}

The variety whose Gromov-Witten invariants we will be interested
in here is $H = H(2,\proj^2)$, the Hilbert scheme of two points in
the plane.  It functorially parametrizes length two subschemes
of $\proj^2$.  So there exists a subscheme $\XX$ in $H \times \proj^2$
whose fiber over a point of $H$ is exactly the
subscheme parametrized by that point.
$$
\begin{CD}
 \XX @>{f}>>  \proj^2\\
@V\nu VV \\
H\\
\end{CD}
$$
In this section, we collect some basic facts about the geometry of $H$.

Given a degree two subscheme, $S \subset \proj^2$, there exists
a unique line containing $S$.  This correspondence gives rise to a 
morphism 
$$\pi : H \rarr \proj^{2*}.$$
The fiber over a point $[L]$ is the set of subschemes contained
in that line.  This is canonically $\Sym^2(L)$.  So we can realize
$H$ as a $\proj^2$ bundle over $\proj^{2*}$.  It is simply
$\proj (\Sym^2(S))$ where $S$ is the tautological rank 2 
subbundle on $\proj^{2*}$.

One thing that we can see immediately from this description of $H$
is that it is a smooth variety.  It follows that $\XX$ is also
smooth, and the morphism from $\XX$ to $H$ is finite, and simply
branched over $\Delta$, the smooth subvariety of $H$ parametrizing
nonreduced length two subschemes.

This representation of $H$ as a projective bundle also gives us
a very good understanding of its cohomology ring.  First of
all, we can conclude 
that the Chow ring and cohomology ring are isomorphic,
so we will use Chow notation from here on to avoid doubling
indices.
We see that $\Pic(H)$ has rank 2 and is generated
by classes $T_1= \pi^*({\rm hyperplane})$ and $T_2 =\oh(1)$.
$T_1$ is the divisor consisting
of schemes whose associated line is incident to a fixed point. 
$T_2$ can be represented by the divisor consisting of all 
subschemes
incident to a fixed line in $\proj^2$.

From these representations, it is clear that both of the classes 
$T_1$ and $T_2$ move without basepoints.
It follows that any effective curve in $H$ has non-negative intersection
with each of them.  
In fact, the converse is true.
If we set $B_1$ to be the curve of subschemes supported at a fixed point,
and $B_2$ to be a line in a fiber of $\pi$,
we see that
$$B_i \cdot T_j = \delta_{ij}.$$
This implies that the cone of effective curves consists exactly
of positive integral combinations of $B_1$ and $B_2$, so we will
fix these as our basis for $A_1H$ and use $(a,b)$ to denote 
$aB_1 + bB_2$ throughout.

We know that $A^*(H)$ is generated by $T_1$
and $T_2$ modulo the relations $T_1^3=0$ and 
$$T_2^3 + c_1 T_2^2 +c_2 T_2=0 $$
where $c_i =\pi^* c_i(\Sym^2(S)).$ (Note that $c_3$ must vanish for
dimension reasons.)  We can easily
compute these Chern classes.  They are
\begin{align*}
c_1& = -3T_1\\ 
c_2& = 6T_1^2.
\end{align*}

If we define $T_3=T_1^2$ and $T_4=T_1T_2-2T_1^2$ then we can extend
these to a Poincare dual basis of the integral Chow ring of $H$ by which we 
mean a basis satisfying
 $$\int T_i \cup T_{8-j} = \delta_{ij}.$$
Here we take $T_0$ to be the fundamental class and $T_8$ to be
the point class.  Our choice of $T_4$
is made because this is the class represented by the cycle
$$\Gamma(p) = \{ {\mbox{subschemes incident to {\it p}}} \} $$
where $p$ is a point of $\proj^2$.  Notice that it is clear that
$\int T_4 \cup T_4=1$ and $\int T_4 \cup T_3 = 0$ as these correspond
to the statements that there is a unique degree two subscheme of $\proj^2$
meeting two general points, and that there are no subschemes contained in a 
general line and meeting a general point.  Also, we have already
constructed $T_6$ and $T_7$ as our basis
for the curve classes.  The only remaining element of the basis
is $T_5 = T_1^2-T_1T_2+T_2^2$ for which we know of no direct geometric
interpretation.  However, the class $S_5 = T_5 + T_3$ is represented by
the closure of the locus $\{ \{p,q\} : p \in l_1, q\in l_2\}$ where
$l_1$ and $l_2$ are distinct lines in $\proj^2$.

The description of $H$ as a projective bundle also allows 
us to compute its canonical
class.  The result is that $c_1(TH)= 3T_2$.  

We will see that it will be especially important for us to understand
how the diagonal, $\Delta$, which parametrizes length 2 subschemes of
$\proj^2$ supported at a single point, sits inside $H$.  
If we think of $H$ as 
$\proj (\Sym^2S)$, then $\Delta$ is the image of 
$\proj (S)$ 
under the degree 2 
Veronese embedding.  From this we can deduce that $\Pic(\Delta)$
is generated by $T_1$ and $\frac{1}{2}T_2$.  Again, these 
are nef divisors so any
curve in $\Delta$ has positive integral intersection with each of them.  
Since $B_1$ and $2B_2$ can both be represented by curves in $\Delta$ we
see that the cone of effective curves in $\Delta$ is just the locus $(a,b)$
with $a$ and $b$ nonnegative, and $b$ even.
By intersecting it
with curves, we can
calculate that as a divisor in $H$, $$\Delta \equiv 2(T_2-T_1).$$

\subsection{Geometry of Hyperelliptic Curves}
\label{ghc}
Although the method we use to solve enumerative problems
is to immediately forget about the actual hyperelliptic
curves and replace them with rational curves in $H$, we will
discuss in this section the natural space of hyperelliptic
maps, and make explicit the relationship between this space
and the genus 0 space on which we calculate.  In addition
to clarifying our strategy, we will eventually make use
of this direct construction to do some of the deformation
theory calculations in Section \ref{deform}.

The fundamental result about hyperelliptic curves which makes
their geometry much easier to study than that of arbitrary curves
is the following.

\begin{lm}
\label{hypvan}
If $f:C \rarr \proj^r$ is a morphism from a hyperelliptic
curve which does not factor through the hyperelliptic
map, then $f^*(\oh(1))$ has no higher cohomology.
\end{lm}
\bpf
In genus 0 or 1 all line bundles of positive degree have no 
higher cohomology, so the statement is vacuous.  In genus 2
or higher, we know that the canonical morphism is 2 to 1 onto
a rational normal curve.  Given a divisor $D$ in the linear
series, Serre duality says that the dimension of $H^1(D)$ is given
by the dimension of $H^0(K-D)$ which is sections of the
canonical bundle vanishing on the divisor, or equivalently,
hyperplanes in $\proj(H^0(K)^{\vee})$ containing the image of
the divisor under the canonical morphism.  
Since any collection of distinct points
on a rational normal curve are in general linear position,
the only way that a divisor can fail
to impose the maximum number of conditions on the canonical
series is for it to contain a pair of hyperelliptically
conjugate points.  But if every element of the linear series
contains a conjugate pair, we can deduce that 
the map must factor through the hyperelliptic map.
\epf

We will make use of this fact to study one natural moduli space
of hyperelliptic plane curves.
This space, which we denote by $H_g(\proj^2,d)$ parametrizes
hyperelliptic curves of genus $g$ mapping to the plane.
It is constructed as a locally closed 
subvariety of $M_{g,0} (\proj^2,d)$.  

Look at the locus $H_g \subset M_g$ consisting of hyperelliptic
curves.  This locus is in fact a smooth substack of $M_g$.  This
follows from the standard realization of
the moduli space of hyperelliptic curves as the quotient by the 
symmetric group of $M_{0,2g+2}$.  In other words, a hyperelliptic 
curve is determined up to the hyperelliptic involution by 
the associated configuration
of branch points on $\proj^1$.  We first define the space of
hyperelliptic maps to be the  
preimage of this hyperelliptic 
locus in $M_{g,0}(\proj^2,d)$.  Call
this space $\tilde H_g(\proj^2,d)$.  
\begin{footnote}
{In genus 0 or 1, the space
of hyperelliptic curves should really be thought of as a smooth
family over the moduli space of curves.  $H_g$ should parametrize the 
choice of curve along with the choice of hyperelliptic involution.
In this case, rather than a preimage,
the hyperelliptic map space should be constructed as the 
fiber product of $H_g$ with the map space.  Throughout
the paper, we will think of hyperelliptic curves in this way.
That is, the term hyperelliptic curve should be taken to mean
a curve together with a choice of hyperelliptic involution.}
\end{footnote}

If we look at the open subset of this space which parametrizes
maps which are birational onto their image, we get a reasonable
candidate for a moduli space of hyperelliptic plane curves.  It
will be convenient for us to look at a slightly larger open subset.
This space, which we will denote by $H_g(\proj^2,d)$, is the subset
parametrizing maps which do not factor through the hyperelliptic 
involution.

In this context, the significance of the 
vanishing result introduced earlier is the following.

\begin{tm} The natural morphism from
$H_g(\proj^2,d)$ to $H_g$ is smooth.  
\end{tm}
\bpf
Consider the Euler sequence:
$$ 0 \rarr \oh \rarr \oh(1)^{\oplus 3} \rarr T\proj^2 \rarr 0.$$
It induces a surjection 
$H^1(C,f^*\oh(1))^{\oplus 3} \rarr H^1(C,f^*T\proj^2)$.
We know by Lemma \ref{hypvan} that 
the first group vanishes, so we can conclude the
smoothness of the morphism by Theorem \ref{relobtm}.
\epf

\begin{cor}
$H_g(\proj^2,d)$ is smooth and irreducible.
\end{cor}

\bpf
The smoothness follows immediately, since $H_g$ itself is smooth.
Since $H_g$ is irreducible, the corollary will follow once
we verify that
the fibers of $\nu$ are 
irreducible.  A fiber is just the set of all maps from a fixed
hyperelliptic curve to $\proj^2$ which do not factor through
the hyperelliptic map.  By associating to a map, the corresponding
line bundle, we see that
this space is fibered over $\Pic^d(C)$.  By Lemma
\ref{hypvan} we can see that its image in $\Pic^d(C)$ is the open
subset consisting of line bundles with vanishing first cohomology.  Over
this locus, the moduli space is just an open subset of 
the bundle of 3-tuples
of sections modulo scalars.  
\epf

We will frequently make implicit use of the 
irreducibility statement 
in that it gives us
a meaningful notion of a generic hyperelliptic plane curve.

\subsection{The Basic Correspondence}
Viewing the space in this way, one would be inclined to
take as a compactification of this space the preimage of
the closure of the hyperelliptic locus.  However, this 
compactification differs markedly from the one we will study.

Instead, we make use of the following correspondence.  Given
a hyperelliptic map, we construct a map from a rational
curve to $H$.  The construction of this map is a simple
application of the universal property of $H$.
A hyperelliptic map gives us the following diagram:

$$
\begin{CD}
 C @>{g}>>  \proj^2\\
@V{2-1}VV \\
\proj^1 \\
\end{CD}
$$

The graph of this correspondence gives a closed subscheme of
$\proj^1 \times \proj^2$ flat over $\proj^1$. (Flatness is
trivial here since $\proj^1$ is a smooth curve, and any
surjective morphism from an irreducible variety to a smooth
curve is flat.)  Provided
that $g$ does not factor through the hyperelliptic map,
the general fiber of the projection of the graph
is a pair of distinct points. Thus the universal property of $H$ gives
us a natural morphism $f:\proj^1 \rarr H$.

Conversely, if we have a map from $\proj^1$ to $H$, pulling back the
universal subscheme gives us a scheme of degree 2 over $\proj^1$, 
together with a map from this scheme to $\proj^2$.  It will not
necessarily be the case that this scheme is what we would usually
think of as a hyperelliptic curve.  For example, if the map takes
$\proj^1$ entirely into the diagonal, this scheme will be everywhere
non-reduced.  If we assume that the original morphism is transverse
to $\Delta$, then this pull-back will be a smooth hyperelliptic 
curve.  This is because $\Delta$ is exactly the branch locus of
the morphism from the universal subscheme to $H$.

The points of intersection of the rational curve in $H$
with $\Delta$ correspond exactly to the branch points
of the associated hyperelliptic curve.
Since we have computed that $\Delta \equiv 2(T_2-T_1)$, 
we can recover the genus
of the hyperelliptic curve 
from the homology class of the corresponding rational
curve via the Hurwitz formula,
$$2g-2= \beta \cdot \Delta = 2(b-a).$$
We can also recover the degree of the hyperelliptic curve by 
intersecting the rational curve with $T_2$.  

We conclude that a generic hyperelliptic plane curve of degree
$d$ and genus $g$ is represented by a rational curve in $H$
of type $(d-g+1,d)$.  To give a more precise statement, we 
restrict our attention to well-behaved curves.  Inside 
$\mooH$ we look
at the open substack $M^{\rm {tr}}_{0,0}(H,\beta)$ 
parametrizing maps from irreducible rational curves which intersect
$\Delta$ transversally.  Correspondingly, in $H_g(\proj^2,d)$ 
we look
at the open subset $H_g^{\rm{tr}}(\proj^2,d)$ which
parametrizes maps from smooth hyperelliptic curves 
such that no two points which are hyperelliptically conjugate map to
the same point, and such that the differential is injective at the branch
points of the hyperelliptic map.  Note that this condition is equivalent
to requiring that the induced map from the hyperelliptic curve
to $\proj^1 \times \proj^2$ is an embedding.

\begin{tm}
\label{annoying}
There is a canonical isomorphism 
$$H_g^{\rm{tr}}(\proj^2,d) \cong M_{0,0}^{\rm{tr}}(H,(d-g-1,d)).$$ 
\end{tm}

\bpf
This follows from a relative version of the construction we
gave above.
Over $H_g^{\rm{tr}}(\proj^2,d)$ we have a smooth family of
hyperelliptic curves, $\CC$.  By considering the morphism induced
by the relative canonical bundle of this family, we get a two
to one map from $\CC$ to a smooth family of rational curves, $\RR$.
$\CC$ also comes with a map to $\proj^2$, and the defining property
of $H_g^{\rm{tr}}(\proj^2,d)$ is that the induced map
from $\CC$ to $\RR \times \proj^2$ is an embedding.
Thus, we have produced over $\RR$ a flat subscheme of $\RR \times \proj^2$.
By the universal property of $H$, we then get a morphism from $\RR$ to $H$.
Now we have produced a diagram
$$
\begin{CD}
\RR @>>>  H\\
@VVV \\
H_g^{\rm{tr}}(\proj^2,d)\\
\end{CD}
$$
which induces the desired morphism from the base to the space
of genus 0 maps to $H$.

Conversely, over $M^{\rm{tr}}_{0,0}(H,(d-g-1,d))$ we have a family
of rational curves with a map to $H$.  Pulling back the universal
subscheme gives us a double cover of this family which maps to $\proj^2$.
The transversality condition means that this cover is branched
over an \'etale multisection of the family, which ensures that it
is a smooth family of hyperelliptic curves.  This gives us the
inverse morphism.

Equivalently, applying the above arguments to families over an arbitrary
base identifies the functors represented by these two moduli spaces.
\epf

\subsection{Curves in the Hilbert Scheme}

At this point, we can be more explicit about our strategy to count
hyperelliptic curves.  By what we have said so far, we can see that
the number of degree $d$ genus $g$ hyperelliptic plane curves passing
through $3d+1$ general points is the same as the number of 
irreducible rational curves of type $(d-g-1,d)$ 
in $H$, transverse to $\Delta $, meeting
$3d+1$ general translates of $\Gamma (p)$.  Hence, we should expect
to find a relationship between our enumerative problem and the
Gromov-Witten invariant $I_{(d-g-1,d) }(T_{4}^{3d+1})$.  What we need to 
do to make this relationship precise is study the structure of
stable maps to $H$ meeting these cycles in order to account for
the contributions of curves which intersect $H$ badly.  
In addition, we would like to ensure that the 
solutions to the enumerative problem are all counted with multiplicity
one.  

To study these questions, we need to consider the
natural action of $\PGL (3)$ on $H$.  Although $H$ is not 
a homogeneous space, it is nearly homogeneous, in that this action has
only 2 orbits -- the dense orbit parametrizing pairs of distinct points,
and the diagonal.  As in the case of homogeneous spaces, one can exploit
this group action in two ways: to get
transversality results, and to control the dimensions and smoothness
of the moduli spaces.

The transversality results will follow from an easy lemma.
We say $X$ is an {\em almost homogeneous space} for a group $G$, if
$X$ is equipped with a $G$ action which has only finitely many orbits. 
These orbits will then give a stratification of $X$.
In this context
there is an analogue of the Kleiman-Bertini Theorem for 
homogeneous spaces.

\begin{lm} 
\label{bigtrans}
Let $X$ be a smooth $G$-almost homogeneous space and
$Y$ a smooth scheme with a morphism $f:Y \rarr X$. Take
$\Gamma$ to be a smooth cycle on $X$ which intersects the orbit 
stratification properly.
Then, for a general $g$ in $G$, 
$f^{-1}(g\Gamma)$ will be pure dimensional of the
expected dimension, and the (possibly empty) open subset 
$f^{-1}(\Gamma_{\rm reg})$ will be smooth, where $\Gamma_{\rm reg}$
is the subset of $\Gamma$ on which the intersection with the orbit 
stratification is transverse.
\end{lm}

\bpf
We know that it is impossible for any component of this fiber product
to have less than the expected dimension.  To see that no component
can have more we can just apply Kleiman-Bertini to each orbit.  Choose
an orbit $O$ of the $G$ action.  Then applying Kleiman-Bertini
to the map $f:f^{-1}(O) \rarr O$ tells us that for a general
translate, this part of the preimage will have at most the
expected dimension.  Here we need the properness of the intersection
to know that the codimension of $\Gamma \cap O$ in $O$ is at least
as large as the codimension of $\Gamma$ in $X$.  Since there are finitely
many orbits, this gives the result.

The proof of the smoothness proceeds in a similar fashion.  On each
orbit we can choose a general translate that ensures that the 
transversality condition is obtained for the restricted map.  Since 
there are only finitely many orbits, we can choose a translate which
is general with respect to each.  Then we know that for each point
of $y$ in $f^{-1}(\Gamma)$, we have
$$D_{f} (T_yY) + T_{f(y)} (\Gamma \cap O) = TO$$
where $O$ is the orbit containing $f(y)$.  However, on $\Gamma_{\rm reg}$
we know in addition that $T_{f(y)}\Gamma + T_{f(y)}O = T_{f(y)}X$
which forces exactly the transversality condition we want.
\epf

Now we turn our attention to understanding how the group action can
be used to control the local structure of the map space.
The expected dimension of the space of unpointed maps of class $\beta$ to
a space $X$ is given by the
formula
$$ d_X(\beta )=-K_H\cdot \beta +\dim(X) - 3.$$
By our calculation of the canonical class of $H$, we conclude
that the expected dimension $d_H(\beta)$ of the moduli spaces $\MonH$
is
$$d_H = 3b + 1.$$

Since $H$ is not a convex space, the moduli spaces $\MonH$ do not 
in general have this expected dimension.  The almost homogeneity
of $H$ gives us unusually good 
control here, though.
The first observation is that if we restrict our attention to
curves which intersect $\Delta$ properly, the map spaces are
as nice as we could possibly hope.

\begin{tm}
\label{ahsmooth}
If $f:C \rarr H$ is a stable map of genus 0 such that no component
of $C$ is mapped entirely into $\Delta$, then at the moduli point
$[f]$, $\MonH$ is smooth of the expected dimension.  
\end{tm}

\bpf
Since the action of $\PGL(3)$ is transitive outside of $\Delta$, we
can conclude that the tangent
bundle, $TH$ is generated by global sections outside of $\Delta$.
Hence, if we have a map
$$f:C \rarr H$$
such that no component of $C$ is mapped entirely into the diagonal, then the 
pullback, $f^*(TH)$, is generically generated by global sections.
On a prestable curve of genus 0, this implies that the higher cohomology
vanishes.  As this forces the obstruction bundle to be trivial, the
result follows immediately.
\epf

Thus, the only possible source
of excess dimension in the map space comes from curves which
have components mapping into the diagonal.
Furthermore, $\Delta$ itself is a homogeneous space, and therefore convex,
so we can calculate exactly the dimension of the space of 
rational curves contained in $\Delta$.  By adjunction we have
 $c_{1}(T \Delta) = 2T_1 + T_2$, so
the dimension of the space of 
rational curves contained in $\Delta$ is given by 
$$d_{\Delta}= 2a + b.$$
Whenever
$d_{\Delta}(\beta)>d_H(\beta)$ we should expect to see moduli spaces of
excess dimension.  We can see from our formulas that $d_{\Delta}(a,b)$ will
be greater than $d_H(a,b)$ 
exactly when $a>b$.  Since $\Delta \equiv 2(T_2-T_1)$
this condition can be rephrased as $\beta \cdot \Delta < 0$. In other words
the excess dimension always comes from components that are trapped in the
diagonal.

It is worth mentioning that this 
phenomenon is not special to the variety $H$. 
Any time we have a
divisor $D$ in a variety $X$ we can write
$$c_1(TX)=c_1(TD) + D,$$
so intersection with $D$ essentially measures the difference in the expected
dimension of curves in $D$ and curves in $X$.

This looks promising, since we are obviously not interested in
curves which have negative intersection with $\Delta$, as they
cannot possibly correspond to the hyperelliptic plane curves that
we ultimately want to count.
Unfortunately, it is not the case that if a curve class $\beta$ has positive
intersection with the diagonal then all components of $\MonH$ have the 
expected dimension.  For example, if we write $\beta$ as a sum of two classes 
$$\beta = \beta_1 + \beta_2$$
with $\beta_2$ representing the class of a curve contained in $\Delta$ 
such that $\beta_2 \cdot \Delta <-2$ and $\beta_1 \cdot \Delta >0$ 
then there is a
component of $\MonH$
whose general element consists of a map 

$$f:C \rarr H$$
where $C=C_1\cup C_2$ is a union of two rational curves meeting in a node
with $f_*[C_1]=\beta_1$ and $f_*[C_2]=\beta_2$.  The dimension of such a 
component will be
$$d_H(\beta_1) + d_{\Delta}(\beta_2) - 2 = 3b_1 +1 + 2a_2 +b_2-2
= 3b-1 - \Delta \cdot \beta_2$$
which is strictly greater than $d_H(\beta)$.

In general, there will be components of the moduli space whose general
element corresponds to a highly reducible curve, and many of these 
components will have dimension as large or larger than the expected 
dimension.  
Any of these components could
conceivably give us undesired contributions to the Gromov-Witten invariants.
However, if we restrict our attention to the invariants arising from
incidence conditions in the plane, almost all of these contributions
are zero.

\begin{tm}
\label{main}
Fix a class $\beta=(a,b)$ in $A_1H$ with $b>0$
and $3b+1$ general points $p_1, \ldots , p_{3b+1}$ in $\proj^2$.
 
(i) There exist at most finitely many irreducible curves of class
$\beta$ which are incident to $\Gamma(p_i)$ for all $i$.

(ii) All such curves will intersect $\Delta$ transversally in points
disjoint from the $\Gamma(p_i)$.

(iii)  Given
an arbitrary stable map in class $\beta$ which is incident to all the
cycles, it contains a unique irreducible component which is not
contained in $\Delta$, and that component is of class $(a_0,b)$ for some
$a_0 \leq a$.  

\end{tm}

Before giving the proof of this statement, we elaborate on its
consequences for stable maps to $H$.  We immediately see that aside from
the distinguished component described above, all other components
are of type $(a',0)$.  As these curves have negative intersection
with $\Delta$, they must be contained in the diagonal.  
If we think of $\Delta$ as a $\proj^1$ bundle
over $\proj^2$ via the map taking a double point to its support,
(different from the projection map used earlier)
this is just $a'$ times the class of a fiber. So we can see that
these curves are all just multiple covers of these fibers.  The fibers
are just the curves described in the first section that consist
of all the schemes supported at a fixed point of the plane.

For this reason, adding a component of type $(a,0)$ to a stable map
can never cause it to be incident to any extra points in the plane.
If the $(a,0)$ component meets $\Gamma(p)$ it must in fact be contained
in $\Gamma(p)$, forcing another component of the curve to meet 
the cycle.  Also, since different fibers are disjoint, all the
fibers must be incident to the distinguished central component.
The source curve looks like a comb, with the $(a_0,b)$ component
as the handle, and the $(a',0)$ components as teeth.

Finally, because each of the $(a',0)$ components is forced to pass 
through one of the finitely many points of intersection of
the distinguished component with the diagonal, there are
only finitely many candidates for the $(1,0)$ curves which
are multiply covered.  In other words, there are only finitely
many potential image curves for stable maps incident to all of
the cycles.

\bpf

We will proceed by induction on the number of components of $C$.  
First we consider the case where $C$ is irreducible.  If $C$ is
not contained in $\Delta$ we know that it moves in the expected dimension,
so (i) and (ii) follow immediately 
from the general position lemma (since the 
$\Gamma(p_i)$ 
intersect $\Delta$ properly, and since a general irreducible curve
will intersect $\Delta$ transversally.)

If, on the other hand, $C$ were contained in $\Delta$, then $C$ would
give rise to a non-reduced subscheme of $\proj^2$ which was supported
on a rational curve of degree $\frac{b}{2}$.  $C$ would meet  
$\Gamma(p_i)$ if and only if this subscheme met $p_i$, 
but a rational curve of degree $\frac{b}{2}$ can meet
at most $\frac{3}{2}b-1$ general points (if $b>0$).

Now, suppose we have a reducible curve which meets the $\Gamma(p_i)$.  
From what we have already said, it follows that it cannot
be contained in $\Delta$.  We also know that it must have
at least one component contained in $\Delta$, since this
is the only possible source of excess dimension.  Hence,
we can write our curve as a union $C_0 \cup C_1$ with the point
of intersection $C_0 \cap C_1$ contained in the diagonal.
By the inductive hypothesis, each $C_i$ can meet at most $3b_i+1$
of the cycles.  Without loss of generality, we can assume that $C_0$ meets
$3b_0+1$ and $C_1$ meets $3b_1$.  Again by induction, we know that
there exist only finitely many image curves for $C_0$, and the
only components of $C_0$ which are contained in the diagonal
are of type $(a',0)$.  This means that all potential points
of $C_0 \cap \Delta$ can be contained in a union of finitely many
cycles of the form $\Gamma(q_i)$.  This would then imply
that $C_1$ would have to meet at least one of the $\Gamma(q_i)\cap \Delta$
and meet an additional $3b_0$ of the cycles.  But this is just asking
it to meet $3b_0+1$ general cycles, and to have one of the points
of intersection in $\Delta$ which is impossible.
\epf

%This result is what gives us a clear interpretation of the
%original enumerative problem in terms of curves in $H$.
%Set $E(d,g)$ to be the number of degree $d$ genus $g$ hyperelliptic
%plane curves passing through $3d+1$ general points.
%One obvious corollary is the following.
%\begin{cor}
%\label{enum}
%If $p_1, \ldots , p_n$ are general points in $\proj^2$, then 
%$E(d,g)$ is equal to the number of irreducible curves of type
%$(d-g-1,d)$ incident to $\Gamma(p_1) , \ldots, \Gamma(p_n)$.
%\end{cor}

%\bpf
%By \ref{} we know that any irreducible curve meeting all of the cycles
%must intersect $\Delta$ properly, and since a general irreducible curve
%which is not contained in $\Delta$ will in fact intersect the diagonal
%transversally, we can conclude that all of these rational curves
%do in fact correspond to degree $d$ genus $g$ hyperelliptic plane curves
%passing through the $p_i$.  Conversely, for general points, no
%hyperelliptic curve passing through all of them will have a node
%such that the two points mapping to it under the normalization are
%hyperelliptic conjugates, or will have a Weierstrass point mapping
%to one of the $p_i$.  This second condition assures us that the
%intersection of the pull-backs of the $\Gamma_i$ to the $n$-pointed
%map space is not contained in the pull back of $\Gamma_i \cap \Delta$
%for any $i$ which ensures that the intersection is reduced by the
%second part of our general position lemma.
%\epf

\section{Virtual Contributions}
\label{deform}
\subsection{Smoothness}
Theorem \ref{main} gives us a very sharp picture of the locus 
$$A= \pi (\rho_1^{-1}\Gamma_1 \cap \cdots \cap \rho_n^{-1}\Gamma_n).$$
Since the only moduli in the choice of a stable map meeting
all of the cycles is the choice of multiple cover of the $(1,0)$
curves,
$A$ is a union
of finitely many components, each of which set theoretically decomposes as
a product 
$$M(a_1)\times M(a_2) \times \cdots \times M(a_n).$$
Here
$M(a)$ is the space parametrizing the data of an $a$-sheeted
cover of $\proj^1$ with a choice of point mapping to 0 (the point
of attachment to the distinguished component.)  This space
is just a 
fiber of the evaluation map $\rho_1:\M_{0,1}(\proj^1,a) \rarr \proj^1$.

Since we really want to count just the irreducible curves, we need to 
be able to determine the contribution that each of these components
makes to the Gromov-Witten invariant.  To do this, we first observe that
we are in the relatively simple 
situation described in Theorem \ref{realgw}.
Namely, we just need to intersect the virtual class with a complete subvariety
which is completely contained in the smooth locus.

We want to prove that in a neighborhood of one of these comb
curves, the moduli space of maps is smooth.  The key step in
proving this is to show
that there are
no first order deformations of such a map which resolve the
nodes.  Once this is proven, the smoothness of the moduli
space follows from the smoothness of the moduli space parametrizing
the handle of the comb, and the smoothness of the spaces
parametrizing the teeth.  It also follows that the locus $A$ defined
above with its natural scheme structure is reduced, by our transversality
lemma.

We write $C = C_0 \cup C_1 \cup \cdots \cup C_n$ where
$C_0$ is the handle, and we assume we are given a map $f:C \rarr H$
which takes $C_0$ to a curve intersecting $\Delta$ transversally, and
the other $C_i$ to multiple covers of $(1,0)$ curves.  We write $p_i$
to denote the node where $C_i$ meets $C_0$.
Intuitively, the reason there can be no smoothings of such a map
is that
if we look at a small 
neighborhood of one of the components in $\Delta$, we see that $C_0$ meets
$\Delta$ only once, while $C_i \cdot \Delta \leq -2$.  A small 
neighborhood of $p_i$ on the smoothed curve would then have
to contribute $-1$ to the intersection with $\Delta$.
This immediately precludes
any actual smoothings, but we need to verify that it also rules out the
existence of a first order smoothing.

Suppose we have a flat 
family, $\pi:\tilde C \rarr S$, over the double point, 
$S=\Spec \C [t]/t^2$,
and a morphism $\tilde f:\tilde C \rarr H$ such that the
restriction to the reduced point is as above.
To say that this family smoothes the node at $p_i$ means that in a 
small neighborhood of $p_i$, $\tilde C$ looks like a first-order 
neighborhood of a planar node.  
This scheme can naturally
be written as a union of a first order neighborhood of each of the 
two components, so $\tilde C$ can naturally be written as
the union of $\tilde C_j$.  By embedding
the first order deformation in a global deformation,
it is easy to see that $\tilde C_i$
is isomorphic to a first order neighborhood of $\proj^1$ in
the total space of $\oh(-1)$. 
Then, 
pulling back $\oh(\Delta)$ along with its tautological section 
gives us a line bundle $L$ on $\tilde C$ whose degree on $C_i$ is 
less than $-1$, but
with a global section, whose restriction to $C_0$ vanishes only 
to order 1 at
$p_i$.  
The intersection of $\tilde C_i$ with $C_0$ is a double point,
so
the
section cannot be zero on $\tilde C_i$.

We are now reduced to
proving that there do not exist any line bundles on $\tilde C_i$ whose
degree on $C_i$ is less than $-1$ and have a non-zero section.  
This will follow from two lemmas.

\begin{lm}
All line bundles on $\tilde C_i$
are pulled back from the projection map to $\proj^1.$  
\end{lm}
\bpf
Twisting up by a high power of $\oh(1)$ we can assume that
our line bundle has a section whose restriction to
$C_i$ has only simple zeroes.  It is clear that any such
Cartier divisor on $\tilde C_i$ extends to a divisor
on the total space of $\oh(-1)$.  Now the result follows
from the familiar fact that the Picard group of the
total space of a vector
bundle is equal to the Picard group of the base.
\epf

We now know that the restriction of $L$ to $\tilde C_i$ is
of the form $\oh(-d)$ for some $d$ greater than or equal to 2.
All that remains is to show that these bundles have no sections.

\begin{lm}
Let $X$ be a scheme, with $E$ and $F$ vector bundles over $X$.
Let $\tilde X$ be the first order neighborhood of $X$
in $E$, with $\pi : \tilde X \rarr X$ the natural projection map.
There is a canonical isomorphism 
$$H^0(\tilde X, \pi^*F) \cong H^0(X,F) \sumo \Hom (E,F).$$
\end{lm}

\bpf
Pullback and restriction of sections give us the sequence
$$H^0(X,F) \rarr H^0(\tilde X, \pi^*F) \rarr H^0(X,F) $$
which shows that $H^0(X,F)$ is a direct summand.  We need
to show that the kernel of the restriction map is $\Hom(E,F)$.
Given an element of $\Hom(E,F)$ we interpret it as a morphism between
the total spaces of the bundles.  Restricting to $\tilde X$ gives
us a section which vanishes on $X$. The inverse of this map is given by
interpreting an element of $H^0(\tilde X, \pi^*F)$ as a morphism over $X$ 
from $\tilde X$ to the total space of $F$.  The differential of this
map gives a morphism of relative tangent bundles from $T_{\tilde X/X}$
to $T_{F/X}$, which are naturally identified with $E$ and $F$
respectively. 

\epf

Applying this lemma with $X=\proj^1$, $E=\oh(-1)$, and 
$F=\oh(-d)$ completes the proof of the smoothness of the moduli space.

\subsection{Tangent-Obstruction Sequence}

In order to apply Theorem \ref{realgw} to compute the contributions
from the positive dimensional components, we need to understand
the obstruction bundle $\EE$ on these loci.  Our only
tool to study this bundle
is the tangent obstruction sequence:
\begin{equation*}
\begin{split}
0 \rarr H^0(TC) \rarr &H^0(f^*TH) \rarr TM \stackrel{\phi}{\rarr} 
\Ext^1(\Omega_C, \oh_C) \rarr \\
&\rarr H^1(f^*TH) \rarr \EE \rarr 0.
\end{split}
\end{equation*}

From what we have already done, we know that the cokernel of $\phi$
has dimension at least $n$, since none of the nodes can be smoothed.
In fact, the cokernel has dimension exactly
$n$.  That is, the family of deformations of the map
surjects onto the topologically trivial deformations of the curve.  This
is equivalent to saying that the deformations of the irreducible curve
realize all possible deformations of 
configurations of $n$ points in $\proj^1$ via 
intersection with $\Delta$.  Translating this back into the language
of hyperelliptic plane curves, it just means that 
the map from the moduli space
of hyperelliptic maps to the space of hyperelliptic curves is
smooth which was proven in Section \ref{ghc}.

This means that the obstruction bundle fits into an exact sequence:
\begin{equation}
\label{exact}
0 \rarr \sumo_1^n L_i \rarr H^1(f^*TH) \rarr \EE \rarr 0.
\end{equation}
Here $L_i$ is the line bundle corresponding to the deformation which
resolves the $i^{\rm th}$ node.

To analyze the middle term of this sequence, we look at the normalization
sequence:
\begin{equation}
\label{normseq}
\begin{split}
0 \rarr &H^0(C,f^*TH) \rarr \sumo_0^n H^0(C_i,f_i^*TH)  \rarr \\
&\rarr \sumo_1^n TH_{p_i} \rarr 
H^1(C, f^*TH) \rarr \sumo_0^n H^1(C_i,f_i^*TH) 
\rarr 0. 
\end{split}
\end{equation}
We know that $H^1(C_0,f_0^*TH) = 0$ by almost homogeneity, and
since $H^1(f_i^*T\Delta) = 0$ for the same reason, we can conclude
that 
$H^1(C_i, f_i^*TH) = H^1(C_i,f_i^* N_{\Delta/H})$.
As $N$ has degree $-2a_i$ on $C_i$, the rank of the last bundle
in sequence \ref{normseq} is $\sum_1^n (2a_i-1)$.

However, since we know by smoothness 
that the rank of $\EE$ is $\sum_1^n (2a_i-2)$, sequence \ref{exact}
forces
$$H^1(C,f^*TH) \cong \sumo_1^n H^1(C_i,f_i^*N).$$  
%Equivalently,
%the argument above showing that the deformations of the map
%surjected onto the topologically trivial deformations of the
%curve also implies the surjectivity of 
%$$H^0(C_0,f_0^*TH) \rarr \sumo_1^n N_{p_i}.$$  As homogeneity
%(or inspection) gives the surjectivity of 
%$H^0(C_i,f_i^*T\Delta) \rarr T_{p_i}\Delta$, we can
%conclude the desired isomorphism.
Hence, we can rewrite the end of the tangent obstruction sequence as
$$0\rarr \sumo L_i \rarr \sumo H^1(C_i,f_i^*N) \rarr \EE \rarr 0.$$
In fact, this sequence naturally splits as a direct sum of
$n$ exact sequences
$$0 \rarr L_i \rarr H^1(C_i,f_i^*N) \rarr \EE_{a_i} \rarr 0$$
where $\EE_a$ is a vector bundle of rank $2a-2$ on $M(a)$.
This bundle is the one that arises when $n=1$.  We won't prove this 
splitting here, since we are only interested in the top Chern
class of $\EE$, and we can already see that the Chern classes of
$\EE$ must be equal to the Chern classes of $\sumo \EE_{a_i}$ since
these bundles fit into the same exact sequence.

\subsection{Evaluation of the Euler Class}
The results of the previous section reduce our problem to
evaluating the top Chern class of $\EE_a$
which is a vector bundle on $M(a)$.  We will see that this Chern 
class vanishes, so
it would be nice to exhibit either a trivial subbundle or trivial
quotient bundle of $\EE_a$.  We have been unable to find such a bundle.
Instead we use a trick.  We will realize $M(a)$ as the moduli space of
maps to another variety, in such a way that the obstruction bundle will
again be $\EE_a$.  The number $c_{\rm top}(\EE_a)$ will then be identified
as a Gromov-Witten invariant of this variety, and we will be able to 
evaluate it using the deformation invariance property.

We consider the variety $X$ obtained by blowing up $\proj^2$ at a point,
and then blowing up a point on the exceptional divisor.  This gives
us a surface with two exceptional divisors $A$ and $B$ meeting in a node.
We take $A$ to be the $-1$ curve and $B$ to be the $-2$ curve  
and set $\beta_a = A + aB$.  All representative stable maps
in class $\beta_a$ consist of a map from a reducible curve with one 
component mapping isomorphically onto $A$ and the rest forming an
$a$ sheeted cover of $B$.  The same arguments as before show that
the moduli space of maps, $\M_{0,0}(X,\beta_a)$ is smooth 
and isomorphic to $M(a)$. 
The expected dimension of these moduli spaces is 0, independent of $a$.
So for each $a$ there is a zero point Gromov-Witten invariant, 
$I^X_{\beta_a}$, which is just the degree of the virtual fundamental
class.  As the whole moduli space is smooth, again
we can realize this class as the top Chern class of a vector bundle
which sits in the exact sequence
$$ 0 \rarr L \rarr H^1(f^*TX) \rarr \EE_a \rarr 0$$
where $L$ is the line bundle corresponding to smoothing the 
unique node lying on the component of the curve which maps to $A$.
Finally, since the normal bundle to $B$ has degree $-2$, we can 
realize the middle term as $H^1(C,f^* \oh(-2))$.  This confirms
that we are indeed looking at the same vector bundles as arise in
determining the Gromov-Witten invariants of $H$ (at least
in K-theory).

Now we are finally in a position to prove the vanishing result we are
after.

\begin{pr} 
\label{vanish}
For all $a\geq2$, $c_{\rm top} (\EE_a) = 0$.
\end{pr}

\bpf
We have already seen that this Chern class is equal to $I^X_{\beta_a}$.
We consider a one parameter family of smooth varieties specializing to
$X$.  We construct it by starting with $\A^1 \times \proj^2$ and first
blowing up $\A^1 \times [1,0,0]$, and then blowing up the proper transform
of the locus $\{(t,[1,0,t])\}$.
If we denote by $X_t$ the fiber of the projection to $\A^1$, 
$X_t$ is the blow-up of $\proj^2$ at $[1,0,0]$ and $[1,0,t]$, and
$X_0 \cong X$.  By deformation invariance, we can
compute our invariant on $X_1$ instead.  We just need to know which
homology class $\beta_a$ corresponds to.  If we label the two exceptional
divisors in $X_t$ as $D_0$ and $D_t$ then $D_0$ specializes to $A+B$
and $D_t$ specializes to $B$.  So the class corresponding to
$\beta_a = A+aB$ is $D_t + a(D_0-D_t)$ which cannot be represented by
any effective curves if $a \geq 2$.  This forces the vanishing of
Gromov-Witten invariants in this class.
\epf

With this result, we can conclude the main theorem
about the relationship between Gromov-Witten invariants of $H$ and 
enumerative geometry of hyperelliptic plane curves.  We set $E(d,g)$
to be the number of hyperelliptic curves passing through $3d+1$ general
points.

\begin{tm}
\label{formula}
The enumerative numbers $E(d,g)$ satisfy the equation
$$ I_{(d-g-1,d)}(T_4^{3d+1}) = \sum_{h\geq g} \binom{2h+2}{h-g}E(d,h). $$
\end{tm}

\bpf
We know that there is a zero dimensional component of the moduli
space corresponding to curves of the following type.  Take a curve
of type $(a',b)$ incident to all of the $\Gamma (p_i)$ and attach
to it $a-a'$ rational curves each mapping isomorphically onto
a $(1,0)$ curve.  The number of such maps is equal to the number
of irreducible $(a',b)$ curves through the cycles times the number
of choices for attachment points of the $(1,0)$ curves.  
The
$(a',b)$ curve will meet the diagonal in $2b-2a'$ points, of which
we choose $a-a'$.  
Also, by part (ii) of Theorem \ref{main} and our general position result,
it follows that these zero dimensional components are reduced, and
so count with multiplicity one.
The formula then follows immediately from the relationship
between $(a,b)$ and $(d,g)$, since we have established that the
only positive dimensional components come from looking at higher degree
multiple covers of the $(1,0)$ curves and these make no contribution
to the Gromov-Witten invariant.
\epf

\subsection{Genus 0 and 1}
\label{g01}

The $E(d,g)$ will all be 0 for curves of genus 0 
or 1, since it is well known that curves of genus 0 can meet at
most $3d-1$ points, and elliptic curves can meet at most $3d$ points.
From our point of view, this is because these curves have extra
$g^1_2$s.  It is still true that the rational curves in $H$ corresponding
to these curves move in a $3d+1$ dimensional family, but the
correspondence between a hyperelliptic curve in $\proj^2$ and a rational
curve in $H$ actually chooses a hyperelliptic curve with a choice
of hyperelliptic involution.  On curves of genus greater than 1, this 
is no choice at all, since such a curve can have at most one 
hyperelliptic involution.  Elliptic curves have a 1-parameter
family of $g^1_2$s, and rational curves have a 2-parameter family.
We can still use the Gromov-Witten invariants of $H$ to compute
the genus 0 and 1 Severi degrees however.  This is done by looking
at a slightly larger set of Gromov-Witten invariants and by considering
the corresponding wider class of enumerative problems.

We will solve the following enumerative problem:  given $k$ general
points, $p_1, \ldots, p_k$, and $l$ general pairs of points,
$q_1,r_1,\ldots, q_l,r_l$, with $k+3l=3d+1$, 
how many hyperelliptic curves of genus $g$ and degree $d$ pass through
all the points, and satisfy the additional condition that for some 
choice of hyperelliptic involution, $q_i$ is hyperelliptically
conjugate to $r_i$ for all $i$. If we set $E^l(d,g)$ to
be the solution to this problem, then we get the following result.
 
\begin{tm}
\label{var}
$$I_{(d-g-1,d)}(T_8^l \cdot T_4^{3(d-l)+1}) = 
\sum_{h\geq g} \binom{2h+2}{h-g} E^l(d,h)$$

\end{tm}
 
\bpf
The connection to Gromov-Witten
invariants is that a hyperelliptic curve will satisfy the
condition that it meets $q$ and $r$ and has the corresponding
points as hyperelliptic conjugates if and only if the
associated rational curve meets the point in $H$ parametrizing
the subscheme $\{q,r\}$. So we can see that the solution
to our enumerative problem is exactly the number of irreducible
curves in $H$ of the appropriate class meeting $l$ general
points, and $k$ general translates of $\Gamma(p)$.  Now
we can apply the same arguments as were used in the proof of
Theorem \ref{formula} to conclude this result.  The only place in our
proof where we needed to refer to the classes that were
being intersected was in Theorem \ref{main} where we showed that 
almost no reducible curves could contribute.  As we choose our
representative of the point class to be a point outside of $\Delta$,
no curves moving in excess dimension can satisfy this condition,
and the exact same argument will hold for the invariants involving
this class. 
\epf

As any pair of points on an
elliptic curve are conjugate under a unique hyperelliptic involution,
and any 2 pairs of points determine a unique hyperelliptic
map on $\proj^1$, we can conclude that the genus 0 and 1 Severi
degrees are given by $E^2(d,0)$ and $E^1(d,1)$ respectively.  
Theorem \ref{var} then gives us a means of computing these
degrees.  In genus 0, of course, one can compute these Severi degrees
much more readily
by working directly with the genus 0 Gromov-Witten invariants of
the plane.  Computing them this way does 
provide a useful check on our calculations, though.

\subsection{Hirzebruch Surfaces}
We remark that the exact same calculations apply to the problem of 
finding the genus 0 Severi degrees of the Hirzebruch surface $\F_2$, 
recovering a result of Abramovich and Bertram, \cite{ab} 
(see \cite{ravi} who computes
Severi degrees of all genera of ruled surfaces, and includes a
discussion of the
result mentioned.)  

We just sketch how this application can be carried out.
The moduli spaces of rational maps to $\F_2$ 
all have the expected dimension, except for
curves which involve a multiple cover of the $-2$ curve.  Just as
in the double blow-up of $\proj^2$, the obstruction theory for
curves with components mapping onto the exceptional divisor is
identical to that of curves in $H$ with components mapping multiply
onto $(1,0)$ curves.  One can conclude that exactly the same
formula as in Theorem \ref{formula} relates the genus 0 Severi degrees to the 
Gromov-Witten invariants on $\F_2$.  In this case, one can 
say even more, since it follows by deformation invariance that the
Gromov-Witten invariants of $\F_2$ are identical to those for
$\F_0$.  As $\F_0 = \proj^1 \times \proj^1$ 
is a homogeneous space, we conclude
that the Gromov-Witten invariants of $\F_2$ are equal to appropriate
genus 0 Severi degrees on 
$\F_0$, so the formula can be interpreted
as relating the genus 0 Severi degrees of the two varieties.  
Set $N_{\F_i}(D)$ to be the genus 0 Severi degree of $\F_i$
for the linear series $D$.  On $\F_i$, take 
$F$ to be the class of a fiber, S to be the class
of a section with self intersection $i$, and $E$ to be the class
of a section with self intersection $-i$.
The result is the formula
$$N_{\F_0}(aS+(b+a)F)=\sum_{i=0}^{a-1}N_{\F_2}(aS+bF-iE).$$

Abramovich and Bertram obtain this result by studying the deformation which
specializes $\F_0$ to $\F_2$.  This is of course closely analogous to the
deformation which we use in the proof of Theorem \ref{vanish}.  In fact,
it was the similarity between their result and Theorem \ref{formula} 
which suggested to us the possibility of a proof via surface geometry.
(Our original proof was along different lines.)

\section{Calculation of Gromov-Witten Invariants}
\label{calc}
\subsection{Two Point Invariants}
\label{2point}
To actually compute the enumerative numbers we are after, we need a method
to calculate the Gromov-Witten invariants of $H$.  The main tool here
will be the First Reconstruction Theorem of \cite{km}.  This 
result (mentioned in Section 1) gives an explicit
algorithm by which all Gromov-Witten invariants of $H$ can be computed
in terms of just the 2 point numbers.  We can apply this theorem 
because we have already seen that the cohomology ring of $H$ is generated
by divisors.  This immediately reduces our problem to finding all
numbers of the form $I_{(a,b)}(\gamma_1 \cdot \gamma_2)$.  Since $H$ is four 
dimensional, each $\gamma_i$ can impose at most 3 conditions on curves,
so together, 2 classes can impose at most 6 conditions.  However, if $b>1$
the expected dimension of curves of type $(a,b)$ is at least 7, so 2 point
numbers exist only for curves of types $(a,0)$ and $(a,1)$.

We start by looking at $(a,0)$ curves.  We have already observed
that all curves of this type are $a$-sheeted covers of curves
of type $(1,0)$. 
The expected dimension of the space of $(a,0)$ curves is 1, independent of 
$a$.  So the only invariants here are of the form $I_{(a,0)}(\gamma)$ where
$\gamma$ is an element of $A^2H$.  

We first compute the invariant on the two geometric codimension two
loci $T_4$ and $S_5$.

\begin{lm}
For all $a$, $I_{(a,0)}(T_4) = I_{(a,0)}(S_5) = 0.$
\end{lm}

\bpf
In each case, the vanishing of the invariant is 
because the class actually imposes 2 conditions on $(a,0)$ 
curves.  Consider the $a=1$ case.  We can see that for a $(1,0)$ curve to
be incident to $\Gamma(p)$, we have to take the curve
to be the curve of nonreduced schemes supported at $p$.  This
condition is obviously codimension 2 in the space of all $(1,0)$ 
curves.  We were only expecting a codimension 1 condition, and as
a result we see that $I_{(1,0)}(T_4)=0$ since $\Gamma(p)$ represents
$T_4$.  Similarly for an $(a,0)$ curve to meet $\Gamma(p)$, the
$(1,0)$ curve which it covers must be the one supported at $p$.
This is again a codimension 2 condition.  

For the $S_5$ invariant, we use
the representation of $S_5$ as the cycle of subschemes incident
to each of 2 lines.  It is easy to see that a $(1,0)$ curve
can meet this cycle only if it is the curve consisting of
subschemes supported at the intersection of the two lines.
Now the same argument that we used for $T_4$ gives us the result.
\epf

This leaves $I_{(a,0)}(T_3)$.  We will compute this directly for $a=1$.
Under the natural identification of $\M_{0,0}(H,(1,0))$ with $\proj^2$,
the locus of curves incident to a general representative
of $T_3$ is a line.  It
follows that $I_{(1,0)}(T_3)$ is simply equal to the degree of the virtual
class. 

In this case, we can directly compute the virtual class.
We have that the moduli space is just $\proj^2$, and the
universal curve over it is the variety of complete flags in $\proj^2$ 
which maps isomorphically
onto $\Delta$.  The virtual class is then 
$c_1(R^1\pi _*(f^* N_{\Delta / H}))$.  As we know that the normal
bundle is strictly negative on all curves in this class, it follows
that $\pi_* f^* N = 0$, so 
$$R^1\pi_* f^*N = -\pi_!(f^*N).$$
We know what the Chern class of $N_{\Delta / H}$ is, since
we have computed its intersection with each of the curve classes.
It is now straightforward to apply Grothendieck-Riemann-Roch to compute
the Chern character of this bundle and read off the desired Chern class.
The result is that $I_{(1,0)}(T_3) = 3$.
We will later calculate $I_{(a,0)}(T_3)$ for $a>1$ 
by means of the associativity
relations.  It is also possible to compute $I_{(1,0)}(T_3)$ by making use 
of more of the associativity relations, thereby avoiding the direct virtual
class calculation.

Now consider the curves of type $(a,1)$.  The main result we want to
prove about these curves is the following:

\begin{tm}
If $a>2$, then all Gromov-Witten invariants for curves of type $(a,1)$
vanish.
\end{tm}

\bpf

We observe that such curves have negative intersection
with the diagonal, but on the other hand, cannot be contained in the diagonal
since 1 is odd.  Hence, all such curves are reducible, and
have a $(0,1)$ curve as one component.  (Irreducible curves of type
$(1,1)$ do not intersect the diagonal, and so cannot be a component
of a connected curve whose other components are all contained in $\Delta$.)
A general $(0,1)$ curve will intersect $\Delta$ in two distinct points,
so a general curve of type $(a,1)$ will be a union of an $(a_1,0)$ curve,
a $(a_2,0)$ curve, and a $(0,1)$ curve where $a_1+a_2=a$.  

We are 
free to choose a basis of cohomology such that every class in our basis
can be represented by cycles which intersect $\Delta$
properly.  Hence, it suffices to
show that any Gromov-Witten invariant involving such classes is
0.  By choosing such cycles, we can guarantee that 
the preimages of general translates will intersect in a subscheme of the
correct codimension, which in this case is 6.  We then
expect the codimension of the image of this intersection in the
unpointed space to be 4.  Of course,
as we saw in our discussion of the $(a,0)$ curves, it is possible
for the image of the intersection to have less than the expected
dimension.  In this case there is nothing to prove, though, since
the invariants are forced to vanish.

We assume then that the image of the intersection of the cycles
has codimension 4 in the space of unpointed maps.
This fact alone is enough to determine very precisely what
the structure of this locus is.  We consider the map
$$g: \M_{0,0}(H,(a,1)) \rarr \M_{0,0}(H,(0,1))$$
which forgets the components contained in the 
diagonal.  Then we observe that if a given point is in the intersection,
that all points in that component of the fiber of $g$ must also be in the
intersection.  This is simply because the difference between stable
maps that lie in the same fiber of $g$ is just the choice of multiple cover
of the associated $(1,0)$ curves.  Choosing a different cover cannot
affect incidence conditions.  Since these fibers are already codimension
4, it follows that the intersection consists entirely of a union of 
finitely many components of fibers of $g$.  

Thus, the structure of the solution sets is of the form 
$$\coprod M(a_1) \times M(a_2)$$
with $a_1+a_2=a$ where these are the same $M(a)$ that occurred in 
the previous section.  The virtual classes here are clearly the same
as they were in that context, so we can conclude that any time
one of the $a_i$ is greater than 1, the contribution from that component
is 0.  If $a>2$, then for some $i$, $a_i>1$, so we are done.  
\epf

To calculate the 2-point numbers for curves of type $(a,1)$ with $a \leq 2$
we simply solve directly the appropriate enumerative problems, being
careful to do our computations with cycles that intersect $\Delta$
properly.  
In the $a=2$ case, we apply the same virtual calculation to disregard
solutions involving double covers of a $(1,0)$ curve.
This is elementary, but somewhat tedious.  We compile the results in a
table.
\\

\begin{tabular}{|r|r|r|r|r|r|r|}
\hline
$I_{(a,1)}(T_i,T_j)$ &$(i,j)=(3,8)$&$(4,8)$&$(5,8)$&$(6,6)$&$(6,7)$&$(7,7)$\\
\hline
a=0 &0&0&1&0&0&1\\
1&1&2&1&1&2&-2\\
2&1&0&0&4&-2&1\\
\hline
\end{tabular}
\\

This leaves only the numbers $I_{(a,0)}(T_3)$, which we will
now compute via the associativity equations.  

We will need to use only the equation associated to the following
diagram

$$\feyn 6312 \ \sim\ \Bigl(
\begin{picture}(40,10)(-20,17)
\put(0,16){\line(0,1){8}}
\put(0,16){\line(3,-1){10}}
\put(0,16){\line(-3,-1){10}}
\put(0,24){\line(3,1){10}}
\put(0,24){\line(-3,1){10}}
\put(12,8){1}
\put(-17,8){3}
\put(12,26){2}
\put(-17,26){6}
\end{picture}\Bigr)$$
with no cohomology classes and $\beta = (a,1)$.

Most of the terms in this equation vanish immediately.
Since $T_2 \cdot (a,0) = 0$,
the divisor axiom forces 
any Gromov-Witten invariant for $(a,0)$ curves containing a 
$T_2$ to vanish.  This already simplifies the recursion to
\begin{equation*}
\begin{split}
I_{(a,1)}(T_3,T_6,T_4+&2T_3)-I_{(a,1)}(T_1,T_3,T_8) - \\
 \sum_{a_1+a_2=a}
&I_{(a_1,0)}(T_1,T_1,T_3) \cdot I_{(a_2,1)}(T_2,T_6,T_7) = 0.
\end{split}
\end{equation*}
Applying this in case $a>2$ we can use the vanishing of $(a,1)$ invariants
and our calculation of $I_{(a_2,1)}(T_6,T_7)$ to reduce to

$$(a-1)^2I_{(a-1,0)}(T_3) = (a-2)^2I_{(a-2,0)}(T_3)$$

Given our computation that $I_{(1,0)}(T_3)=3$, this 
inductively determines that $I_{(a,0)}(T_3)=3/a^2$.

\subsection{Enumerative Results}
The results of the previous section give us 
complete knowledge of all 2-point numbers on $H$, so
we can recursively determine arbitrary Gromov-Witten invariants by means 
of the First Reconstruction Theorem.
This algorithm was implemented on
Maple.  We list here some of the results and 
their enumerative consequences.
\\

\begin{tabular}{|r|r|r|r|r|r|r|}
\hline
$I(T_4^{3d+1})$ & d=2 & 3 & 4 & 5 & 6 & 7 \\
\hline
g=0 & 0 & 0 & 405 & 560385 & 1096808499 & 3292618732704 \\
1 & * & 0 & 162 & 224910 & 460743174 & 1470159619803\\
2 & * & * & 27 & 37935 & 89898984 &338090337018\\
3 & * & * & * & 135 & 3933549 & 29267016849\\
4 & * & * & * & * & 405 & 539160678\\
5 & * & * & * & * & *  & 945 \\
\hline
%\end{tabular}

%\begin{tabular}{|r|r|r|r|r|r|r|}
\hline
$E(d,g)$ & & &  &  &  &  \\
\hline
g=0 & 0 & 0 & 0 & 0 & 0 & 0\\
1 & * & 0 & 0 & 0 & 0 & 0\\
2 & * & * & 27 & 36855 & 58444767 & 122824720116 \\
3 & * & * & * & 135 & 3929499 & 23875461099\\
4 & * & * & * & * & 405 & 539149338\\
5 & * & * & * & * & * & 945 \\
\hline
\end{tabular}
\\

The zeros occur in the genus 0 and 1 rows of the
enumerative table because of the extra hyperelliptic
$g^1_2$'s
on rational and elliptic curves.
Since every genus 2 curve is hyperelliptic,
 $E(d,2)$ is the genus 2 Severi degree in the
linear series of degree $d$ plane curves.  These numbers
agree with the calculations of \cite{ch} and have
now been computed again by a different method in \cite{rp}.
$E(4,2)=27$ is the degree of the quartic
discriminant hypersurface. $E(5,2)=36855$ can also be computed
via the 4-nodal formula due to I. Vainsencher.
\\

\begin{tabular}{|r|r|r|r|r|r|r|}
\hline
$I(T_4^{3d-2}*T_8)$ & d=2 & 3 & 4 & 5 & 6 & 7 \\
\hline
g=0 & 0 & 4 & 975 & 500070 & 510209009 & 936943088028\\
1 & * & 1 & 255 & 147780 & 172751014 & 358483479813\\
2 & * & * & 5 & 10138 & 21081609 & 61683241918\\
3 & * & * & * & 12 & 558749 & 3685184208\\
4 & * & * & * & * &  22 & 32184102\\
5 & * & * & * & * & * & 35 \\
\hline
%\end{tabular}

%\begin{tabular}{|r|r|r|r|r|r|r|}
\hline
$E^1(d,g)$ &  &  &  &  &  &  \\
\hline
g=0 & 0 & 0 & 0 & 0 &  0 & 0\\
1 & * & 1 & 225 & 87192 & 57435240 & 60478511040\\
2 & * & * & 5 & 10042 & 16612387 & 33328207904 \\
3 & * & * & * & 12 & 558529 & 3363345078\\
4 & * & * & * & * & 22 & 32183682\\
5 & * & * & * & * & * & 35 \\
\hline
\end{tabular}
\\

The genus 0 row still vanishes. 
As we mentioned in section \ref{g01}, the row $E^1(d,1)$ 
gives the genus 1 Severi degrees of $\proj^2$.
These elliptic numbers agree with computations by
E. Getzler in \cite{ezra} (who also checked $E_{1,6}=57435240$
with the algorithm of L. Caporaso and J. Harris, \cite{ch}).
\\

\begin{tabular}{|r|r|r|r|r|r|r|}
\hline
$I(T_4^{3d-5}*T_8^2)$ & d=2 & 3 & 4 & 5 & 6 & 7 \\
\hline
g=0 & 1 & 16 & 1279 & 317408 &  187613888 & 222541278466\\
1 & * & 1 & 167 & 63228 & 49635964 & 72095337199\\
2 & * & * & 1 & 2536 & 4254399 & 9650092804\\
3 & * & * & * & 1 & 65417 & 402592233\\
4 & * & * & * & * & 1 & 1900762\\
5 & * & * & * & * & * & 1\\
\hline
%\end{tabular}

%\begin{tabular}{|r|r|r|r|r|r|r|}
\hline
$E^2(d,g)$ &  &  &  &  &  &  \\
\hline
g=0 & 1 & 12 & 620 & 87304 &  26312976 & 14616808192\\
1 & * & 1 & 161 & 48032 & 25417860 & 22151587040\\
2 & * & * & 1 & 2528 & 3731098 & 6495881498\\
3 & * & * & * & 1 & 65407 & 383584667\\
4 & * & * & * & * & 1 & 1900750\\
5 & * & * & * & * & * & 1 \\
\hline
\end{tabular}
\\

$E^2(d,0)$ equals $N_d$, the number of degree $d$ rational
curves passing through $3d-1$ general points in $\bold{P}^2$.
The numbers $N_d$ agree with the computation
of M. Kontsevich.

\subsection{Small Quantum Cohomology}

In this section we use our results to give an explicit
presentation of the small quantum cohomology ring of $H$.
We will see that essentially all of the work has been
done in the computations of Section \ref{2point}.  The
general principle at work is that if a variety has
its cohomology generated by divisors, then the
small quantum cohomology can be described in terms
of the two point Gromov-Witten invariants.  From
one point of view, this is an immediate consequence
of the reconstruction theorem used earlier.  However,
in that context there were infinitely many invariants
to be determined, and we merely had a recursive procedure
to find them.  Here, because we know that we just need
to find a finite set of relations, the problem can
be solved in closed form.

We remind the reader of the definition of the small quantum product.  
We introduce variables
$q_1$ and $q_2$, and define a multiplication
on $A^*(H)[[q_1,q_2]]$.  Given $\gamma_1$ and 
$\gamma_2$ in $A^*(H)$ we set
$$\gamma_1 * \gamma_2 = 
\sum I_{(a,b)}(\gamma_1\cdot\gamma_2\cdot T_i)q_1^aq_2^bT_{8-i}.$$
The product is extended to the whole power series by linearity
over $\Q[[q_1,q_2]].$  The associativity of this product
is a consequence of the WDVV equations.

First we compute the quantum product of divisors with
elements of $A^1$ and $A^2$.  (We can get by without
computing other products because the relations in
the ordinary cohomology ring are of degree 3.) 
The result is as follows.
\begin{eqnarray*}
T_1 * T_1 &=& (1-3f)T_3 + 3fT_5\\
T_1 * T_2 &=& 2T_3 + T_4\\
T_2 * T_2 &=& T_3 + T_4 + T_5\\
T_1 * T_3 &=& 3fT_7 + q_1 q_2 + 2q_1^2 q_2\\
T_1 * T_4 &=& T_6 + 2q_1 q_2\\
T_1 * T_5 &=& 2T_6 + (1-3f)T_7 + q_1 q_2\\
T_2 * T_3 &=& T_6 + q_1 q_2 + q_1^2 q_2\\
T_2 * T_4 &=& T_6 + T_7 + 2q_1 q_2\\
T_2 * T_5 &=& T_6 + 2T_7 + q_2 + q_1 q_2
\end{eqnarray*}
where
$$f=\frac {q_1}{1-q_1}= q_1+q_1^2+q_1^3+ \cdots$$
This table follows immediately from the definition of
quantum product, the analogous calculation in the
ordinary cohomology ring, the divisor axiom, and the
two point invariants computed in Section \ref{2point}.
Note that unlike the situation with Fano manifolds,
the quantum multiplication is not defined at the
polynomial level, but only in terms of formal power 
series.

From this we can explicitly write all triple products of divisors
in terms of $T_6$ and $T_7$.  It is then just some linear
algebra to deduce the following relations.
$$T_1*T_1*T_1 = 9f^2T_1*T_2*T_2 - (9f^2-2f)T_2*T_2*T_2 +q_1q_2(q_1-1)$$
$$(1-18f)T_2*T_2*T_2 -3(1-6f)T_1*T_2*T_2 + 6T_1*T_1*T_2 = q_2(q_1^2-2q_1+1)$$
As these are deformations of the two relations defining the 
ordinary cohomology ring of $H$, it follows that they give
a complete set of relations for the quantum ring.  It is
interesting that by clearing denominators, it is in fact possible
to write these relations as {\em polynomials} and not just
power series.  However,
the $T_i$ would not be in the polynomial ring generated by 
the divisors, only in the ring localized at $1-q_1$.

%\bibliographystyle{amsplain}
%\bibliography{database}

\begin{thebibliography}{[FP]}

\bibitem[AB]{ab} D. Abramovich and A. Bertram, {\it private communication}.

\bibitem[B]{b} K. Behrend, {\em Gromov-Witten invariants in algebraic
geometry}, Invent. Math. {\bf 127}, (1997), 601-617.

\bibitem[BF]{bf} K. Behrend and B. Fantechi, 
{\em The intrinsic normal cone}, Invent. Math. {\bf 128}, 
(1997), 45-88.

\bibitem[BP]{rp} P. Belorousski and R. Pandharipande, 
{\em A descendant relation in genus 2}, Preprint, alg-geom/9803072.

\bibitem[CH]{ch} L. Caporaso and J. Harris, {\em Counting
plane curves of any genus}, Invent. Math. {\bf 131}, (1998), 345-392.

\bibitem[FP]{fp} W. Fulton, R. Pandharipande, {\em Notes on
stable maps and quantum cohomology}, in {\em Algebraic Geometry---Santa Cruz
1995}, Amer. Math. Soc., (1997), 45-96.

\bibitem[G]{ezra} E. Getzler, {\em Intersection theory of $\M_{1,4}$
and elliptic Gromov-Witten invariants}, J. Amer. Math. Soc., {\bf 10},
(1997), 973-998.

\bibitem[K]{kollar} J. Koll\'ar, {\em Rational Curves on Algebraic
Varieties}, Springer-Verlag, (1996).

\bibitem[KM]{km} M. Kontsevich and Yu. Manin, {\em Gromov-Witten
classes, quantum cohomology, and enumerative geometry}, Commun. Math.
Phys.  {\bf 164} (1994), 525-562.

\bibitem[LT]{lt} J. Li and G. Tian, {\em Virtual moduli cycles and
Gromov-Witten invariants of algebraic varieties}, J. Amer.
Math. Soc., {\bf 11}, (1998), 119-174.

\bibitem[V]{ravi} R. Vakil, {\em Counting curves of any genus on rational
ruled surfaces}, Preprint, alg-geom/9709003.

\end{thebibliography}
%********

\end{document}